\documentclass[12pt]{article}
\usepackage{amsmath}
\usepackage{latexsym,amssymb}
\usepackage{amsmath,amsthm}
\usepackage{pstricks}
\usepackage{delarray}
\usepackage{tabularx}
\usepackage{amsbsy}
\usepackage{graphicx}
\usepackage{amsfonts}
\usepackage{layout}
\usepackage{mathrsfs}

\input cyracc.def

\newcommand{\per}{\operatorname{per}}
\newcommand{\adj}{\operatorname{adj}}

\begin{document}


\begin{center}
{\Large \bf On permanents of matrices over a commutative additively
idempotent semiring }
\end{center}
\begin{center}
{ \renewcommand\baselinestretch{2.5}
Yan Huang\hspace{0.5cm} Haifeng Lian* \\
{\scriptsize Department of mathematics, Fujian
Agriculture and Forestry University, Fuzhou, 350002, China;\\
Email: xiaowuyav@126.com,   lianhaif@qq.com
 } }
\end{center}
\noindent{\bf Abstract}: Let $R$ be a commutative additively
idempotent semiring. In this paper, some properties and
characterizations for permanents of matrices over $R$ are
established, and several inequalities for permanents are given.
Also, the adjiont matrices of matriecs over $R$ are considered.
Partial results obtained in this paper generalize the corresponding
ones on fuzzy matrices, on lattice matrices and on incline matrices.

\vspace{3mm} {\noindent\bf{Keywords}:}  Permanent; adjoint matrix;
additively idempotent semiring.

\vspace{3mm} \noindent{\bf 2010 Mathematics Subject
Classification}\quad 15A15, 15B15, 16Y60

{\section{{ \bf Introduction}}}

\indent A {\it semiring}  is an algebraic system $(R,+,\cdot)$ in
which $(R,+)$ is an abelian monoid with identity element 0 and
$(R,\cdot)$ is another monoid with identity element $1(\neq 0)$. In
addition, operations $+$ and $\cdot$ are connected by distributivity
and 0 annihilates $R$. A semiring is commutative if $ab=ba$ for all
$a,b\in R$. A semiring $R$ is called an \emph{incline} if $a+1=1$
for all $a\in R$.  Every Boolean algebra, the fuzzy algebra
${\mathbb{F}}=([0,1],\vee,T)$, where $\vee=\max$ and $T$ is a t-norm
(for t-norm, see \cite{KM}), and any  bounded distributive lattice
are examples of inclines.

 A semiring  $R$ is called \emph{an additively
idempotent semiring} if $a+a=a$ for all $a\in R$. Clearly, an
incline is
 additively idempotent.  Additively idempotent semirings are useful
tools in diverse areas such as fuzzy set theory  and decision
analysis; data analysis and preference modeling; classical and
non-classical path-finding problems in graphs; analysis and control
of discrete-event systems (see \cite{GM}).

Permanent of a matrix made its first appearance in the famous
memoirs of Binet\cite{B}. Since then, a large number of works on
permanent theory have been published (see \cite{G,K,LB,M,S,Z}). The
\emph{adjoin matrix} of $A\in M_{n}(R)$ is denoted by $\adj(A)$,
which is defined as an $n\times n$ matrix whose $(i,j)$th entry is
the permanent of $A(j|i)$,
 where $A(j|i)\in M_{n-1}(R)$
obtained from A by deleting the \emph{j}th row and the \emph{i}th
column. Since the late 1980s, many authors have studied adjoint
matrix over
 special cases of additively idempotent semiring. For example, Han and Li
 studied the
 properties of the adjoint matrix of incline matrices and present
Crammer's rule \cite{H1}. Kim et al. \cite{K} studied permanent
theory for fuzzy square matrices and proved that $\per(A\adj
(A))=\per(A)=\per(\adj(A)A)$. This result was generalized to
$D_{01}$-lattice matrices by Zhang  \cite{Z}.
 In 2004, Duan \cite{Du} studied permanents and adjoint matrices of incline
  matrices and proved that $A^n$ is equal
to the adjoint matrix of $A$ if the matrix $A$ satisfies $A\geq
I_n$ and posed the following open problems:\\
{\bf\noindent Problem 1.1.}  Does  the equality
$$\per(A\adj( A))=\per((\adj (A))A)=\per (A)$$
 hold for any square matrix $A$ over an incline?\\
{\bf\noindent Problem 1.2.} Does the equality
$$A^{n-1} =\adj (A)$$ hold for an $n\times n$  matrix $A$ over an
incline satisfying $a_{ii}\geq a_{jk}$ for all $i, j, k \in
\underline{n}$?

In this paper, we will discuss permanents of matrices over a
commutative additively idempotent semiring. In section 2, we give
some properties and characterizations for permanents and obtain some
inequalities for the permanents. In section 3, we consider the
adjoint matrices of matrices over commutative additively idempotent
semiring and answer Problem 1.1 and 1.2. Partial results obtained in
this paper generalize the corresponding results on fuzzy matrices in
\cite{K,RE}on fuzzy matrices, on lattice matrices in \cite{Z} and on
incline matrices in \cite{Du}. Finally, in section 4, we prove
theorem 3.5 stated in section 3.

\vspace{3mm}{\section{{\bf  Properties and characterizations for
permanents}}} \setcounter{section}{2} \setcounter{equation}{0}

\indent In this paper, we always assume that $R$ is a commutative
additively semiring and $\leq$ be the canonical preorder realation
induced by the $+$ operation ($a\leq b$ if and only if $a+b=b$).
Clearly, $0$ is the least element in $R$, i.e., $0\leq a$ for all
$a\in R$. Moreover, for any $a, b, c, d\in R$,
$$\begin{array}{ll} (\textrm{P1})& a\leq b\textrm{ and } c\leq d\quad  \Rightarrow\quad
a+c\leq b+d \textrm{ and } ac\leq bd;\\
(\textrm{P2})& a\leq b\textrm{ and } b\leq a\quad  \Rightarrow\quad
a=b.\end{array}$$ Let $M_{n}(R)$ be the set of all $n\times n$
matrices over $R$. Matrix operations on $M_{n}(R)$ are defined the
same as that in a field. It is easy to verify that $(M_n(R), +,
\cdot)$ is an additively idempotent semiring, thus the canonical
preorder realation $\leq$ ($A\leq B$  if and only if  $a_{ij}\leq
b_{ij}$ for all $i,j\in \underline{n}$) on $M_{n}(R)$ satisfies (P1)
and (P2). Denote by $A^T$ the {\itshape transpose} of $A$.

Denote $\underline{n}=\{1,2,\cdot\cdot\cdot,n\}$. For $A\in M_n(R)$,
the {\itshape permanent} of $A$ is defined as
$$\per(A)=\sum\limits_{\sigma \in S_n}
a_{1\sigma(1)}a_{2\sigma(2)}\dots a_{n\sigma(n)},$$  where $S_n$
denotes the symmetric group of the set $\underline{n}$.
 Checking direct by the definition, we have the following proposition.
\vspace{3mm}

 \noindent {\bf Proposition 2.1}\quad {\itshape Let $A,
B, C\in M_n(R)$ and $\lambda \in R$. Then

(1) $\per(\lambda A)=\lambda^n\per(A)$,  where $\lambda A=(\lambda
a_{ij})_{n\times n}$;

(2) $\per(A^T)=\per(A)$;

(3) $\per(PAQ)=\per(A)$, where $P$ and $Q$ are $n\times n$
permutation matrices;

(4) if $a_{ii}\geq a_{ik}$ for all $i, k \in \underline{n}$, then
$\per(A)=a_{11}a_{22}\cdots a_{nn}$;

(5) $\per  \left(
\begin{array}{cc}
A& C\\
O& B
\end{array}
\right)$ $=\per(A)\per(B)$, where $O$ denotes the $n\times n$ zero
matrix;

(6)  if $A\leq B$, then  $\per(A)\leq \per(B)$.} \hfill $\Box$

\vspace{3mm} For $k\in \underline{n-1}$, set $\Omega_{k,
n}=\{\omega=(i_1, \ldots , i_k )\mid 1< i_1 < \cdots <i_k <n \}$.
For $\alpha=(i_1, \ldots , i_k ) , \beta=(j_1,  \ldots , j_k )\in
\Omega_{k,n}$ and $A\in M_n(R)$, we denote by $A[\alpha|\beta]$ the
$k\times k$ submatrix of $A$ whose $(u,v)-$entry is equal to
$a_{i_uj_v}$ and by $A(\alpha|\beta)$ the $(n-k)\times (n-k)$
submatrix of $A$ obtained from $A$  by deleting rows $\alpha$ and
columns $\beta$. Similar with the Laplace's theorem over a field, we
have the following proposition.

\vspace{3mm} {\bf\noindent Proposition 2.2}\quad {\itshape For $A\in
M_n(R)$ and $\alpha \in \Omega_{k,n}$, we have
$$\per(A)=\sum_{\beta\in \Omega_{k,n}}\per(A[\alpha|\beta
])\per(A(\alpha |\beta)).\eqno(2.1)$$ In particularly, for  $i\in
\underline{n}$, we have
$$\per(A)=\sum_{j=1}^na_{ij}\per(A(i|j)).\eqno(2.2)$$}
\noindent {\bf Proof.} Let $\alpha =(\alpha_1,\alpha_2, \ldots ,
\alpha_k)$. For any $\beta\in \Omega_{k,n}$, we have
 $$\per(A[\alpha| \beta])=\sum_{\sigma\in S_{\alpha,\beta}}\prod^{k}_{i=1}a_{\alpha_i\sigma(\alpha_i)},$$
where $S_{\alpha,\beta}$ is the set of all bijective mappings from
$\alpha $ to $\beta$. Then, there are in total $k!$ terms in
 the expansion of $\per(A[\alpha|\beta])$. Similarly, there are in total $(n-k)!$ terms
in the expansion of $\per(A(\alpha|\beta))$. Since  the product of
any term
 of $\per(A[\alpha|\beta])$ and any term of $\per(A(\alpha|\beta))$ is a term of $\per(A)$,
 there are in total  $k!(n-k)!$ distinct terms in
$\per(A[\alpha|\beta])\per(A(\alpha|\beta))$. Thus, there are in
total $\binom{n}{k}k!(n-k)! =n!$ distinct terms in
$\sum\limits_{\beta\in \Omega_{k,n}}\per(A[\alpha|\beta
])\per(A(\alpha |\beta))$. Since any term of $\sum\limits_{\beta\in
\Omega_{k,n}}\per(A[\alpha|\beta ])\per(A(\alpha |\beta))$ is a term
of $\per(A)$ and $\per(A)$ has $n!$ terms, we have
$$\per(A)=\sum\limits_{\beta\in \Omega_{k,n}}\per(A[\alpha|\beta
])\per(A(\alpha |\beta)),$$ as required.  Taking $\alpha=(i)$ in
(2.1),
 then we have $$\per(A)=\sum_{j=1}^na_{ij}\per(A(i|j)),$$
as required. This completes the proof. \hfill $\Box$ \vspace{3mm}

\noindent {\textbf{Proposition 2.3}}\quad {\itshape For $A,B\in
M_n(R)$, we have
$$\per(AB)\geq \per(A)\per(B).\eqno(2.3)$$}
\noindent {\bf Proof.} Let $T_n$ denote the set of all mappings from
the set $\underline{n}$ to itself. By the definition of permanent,
we have $$\begin{array}{ll}\vspace{2mm}
\per(AB)&=\sum\limits_{\sigma\in
S_n}(AB)_{1\sigma(1)}(AB)_{2\sigma(2)}\dots (AB)_{n\sigma(n)}\\
\vspace{2mm}& =\sum\limits_{\sigma\in S_n}(\sum\limits_{1\leq k\leq
n}a_{1k}b_{k\sigma(1)}) (\sum\limits_{1\leq k\leq
n}a_{2k}b_{k\sigma(2)})\dots (\sum\limits_{1\leq k\leq
n}a_{nk}b_{k\sigma(n)})\\ \vspace{2mm}&=\sum\limits_{\sigma\in
S_n}\sum\limits_{1\leq k_1,k_2,\ldots , k_n\leq
n}a_{1k_1}a_{2k_2}\cdots
 a_{nk_n}b_{k_1\sigma(1)}b_{k_2\sigma(2)}\dots b_{k_n\sigma(n)}\\ \vspace{2mm}
 &=\sum\limits_{\tau \in T_n}(a_{1\tau (1)}a_{2\tau(2)}\cdots
a_{n\tau(n)}\sum\limits_{\sigma\in
S_n}b_{\tau(1)\sigma(1)}b_{\tau(2)\sigma(2)}\dots
b_{\tau(n)\sigma(n)})\\ \vspace{2mm}&=\sum\limits_{\tau \in
S_n}(a_{1\tau (1)}a_{2\tau(2)}\cdots
a_{n\tau(n)}\sum\limits_{\sigma\in
S_n}b_{\tau(1)\sigma(1)}b_{\tau(2)\sigma(2)}\dots
b_{\tau(n)\sigma(n)})\\ \vspace{2mm}&\quad +\sum\limits_{\tau \in
T_n-S_n}(a_{1\tau (1)}a_{2\tau(2)}\cdots
a_{n\tau(n)}\sum\limits_{\sigma\in
S_n}b_{\tau(1)\sigma(1)}b_{\tau(2)\sigma(2)}\dots
b_{\tau(n)\sigma(n)})\\ \vspace{2mm}&\geq \sum\limits_{\tau \in
S_n}(a_{1\tau (1)}a_{2\tau(2)}\cdots
a_{n\tau(n)}\sum\limits_{\sigma\in
S_n}b_{\tau(1)\sigma(1)}b_{\tau(2)\sigma(2)}\dots
b_{\tau(n)\sigma(n)})\\
\vspace{2mm}&= \sum\limits_{\tau \in
S_n}(a_{1\tau (1)}a_{2\tau (2)}\cdots a_{n\tau (n)}\per(B))\\
\vspace{2mm}&=( \sum\limits_{\tau \in S_n}a_{1\tau (1)}a_{2\tau
(2)}\cdots a_{n\tau (n)})\per(B)\\
&=\per(A)\per(B),\end{array}$$ as required. This completes the
proof.  \hfill $\Box$ \vspace{3mm}

\noindent {\bf Remark 2.4}\quad  {\itshape The inequality (2.3) is
not always true on equal. For example, let
$R=([0,+\infty],\vee,\cdot)$, where $\vee=max$ and $\cdot$ is the
usual multiplication of real numbers. It is easy to verify that $R$
is a commutative additively idempotent semiring. Choose
$$
A=\left(\begin{array}{ll}
1&0.5\\
2&2
\end{array}\right),\
B=\left(\begin{array}{ll}
2&1\\
1&0
\end{array}\right),$$
we have $\per(A)=2,\per(B)=1,\per(AB)=4$, so $\per(AB)>
\per(A)\per(B)$. } \vspace{3mm}

\noindent {\bf Corollary 2.5}\quad  {\itshape Let $A\in M_n(R)$. If
$A$ is idempotent (i.e. $A^{2}=A$) with $\per(A)\geq 1$, then $(\per
(A))$ is idempotent too.}

\noindent {\bf Proof.} Since $A^2=A$, using (2.3), we have $(\per
(A))^{2}\leq \per(A)$. Since $\per(A)\geq 1$, we have $(\per
(A))^{2}\geq \per(A)\cdot 1= \per(A)$. Thus  we have $(\per
(A))^{2}=\per(A)$, as required. This completes the proof. \hfill
$\Box$\vspace{3mm}

\vspace{3mm}{\section{{\bf Adjoint matrix of a square matrix }}}
\setcounter{section}{3} \setcounter{equation}{0}

\indent In this section, we will discuss the adjoint matrix of a
square matrix over   a commutative additively idempotent semiring
$R$ and answer Problems 1.1 and 1.2. Recall that the \emph{adjoin
matrix} of $A\in M_{n}(R)$ is denoted by $\adj(A)$, which is defined
as an $n\times n$ matrix whose $(i,j)$th entry is $\per(A(j|i))$,
 where $A(j|i)\in M_{n-1}(R)$
obtained from A by deleting the \emph{j}th row and the \emph{i}th
column. \vspace{3mm}

\noindent  {\textbf{Proposition 3.1}} \quad {\itshape For $A,B\in
M_n(R)$, we have \vspace{2mm}

\indent (1)  $A\leq B$ implies  $\adj(A)\leq \adj(B)$;  \vspace{2mm}

\indent (2)  $\adj(A)+\adj(B)\leq \adj(A+B)$;  \vspace{2mm}

\indent (3)  $(\adj(A))^T=\adj(A^T)$.} \vspace{3mm}

\noindent {\bf Proof.} (1) It follows from the definition of the
adjoin matrix and  proposition 2.1(6).

  (2) Since $A\leq A+B$ and $B\leq A+B$, using (1), we have $\adj(A)\leq
\adj(A+B)$ and $\adj(B)\leq \adj(A+B)$, which implies  $\adj(A)+\adj(B)\leq \adj(A+B)$.

 (3) It is obvious. \hfill $\Box$
\vspace{3mm}

We say  $A\in M_n(R)$ with $n\geq2$ satisfies condition $(*)$ for
convenient, if
$$a_{ii}\geq a_{jk} \textrm{ for all } i,j,k\in \underline{n}.$$
Clearly, if $A$ satisfies condition $(*)$, then so is $A^l$ for any
positive integer $l$.

\vspace{3mm}
 \noindent {\bf Lemma 3.2}\quad  {\itshape Let $A\in
M_n(R)$ with $n\geq2$. If $A$ satisfies condition $(*)$, then for
positive integer $l$, we have
$$\per(A^{l})=(\per A)^l.$$}
\noindent {\bf Proof.}   Since $a_{ii}\geq a_{jk}$ for all $i,j,k\in
\underline{n}$, we have $a_{11}=a_{22}=\cdots=a_{nn}\geq a_{ij}$ for
$i,j\in \underline{n}$. Thus, we have $\per(A)=\sum\limits_{\sigma
\in S_n} a_{1\sigma(1)}a_{2\sigma(2)}\dots a_{n\sigma(n)}=a_{11}^n$.
Denote by $a^{(l)}_{ij}$ the $(i,j)$th entry of $A^l$ with $l\geq
2$, we have
$$a^{(l)}_{ij}=\sum_{j_1,j_2,\cdots,
j_{l-1}\in\underline{n}}a_{ij_1}a_{j_1j_2}\cdots a_{j_{l-1}j}\leq
a_{11}^l $$ and $$a^{(l)}_{ii}=\sum_{j_1,j_2,\cdots,
j_{l-1}\in\underline{n}}a_{ij_1}a_{j_1j_2}\cdots a_{j_{l-1}i}=
a_{11}^l, $$ for $i,j\in\underline{n}$. Thus, we have
$\per(A^l)=\sum\limits_{\sigma \in S_n}(\
a_{1\sigma(1)}^{(l)}a_{2\sigma(2)}^{(l)}\dots
a_{n\sigma(n)}^{(l)})=a_{11}^{ln}=(\per A)^{l}$, as required. This
completes the proof. \hfill $\Box$\vspace{3mm}

 For a max-min $n\times n$ fuzzy matrix $A$, Thomason
\cite{Th} proved that $$\adj(A)=A^{n-1}\eqno(3.1)$$ if $A$ satisfies
condition $(*)$. Duan
 \cite{Du} proved that (3.1) holds for  matrices over distributive lattice with same condition and  matrices over
 commutative incline under stronger condition.
For matrices  over a commutative additively idempotent  semiring
$R$, we have the following theorem. \vspace{3mm}

\noindent {\bf Theorem 3.3} \quad {\itshape Let $A\in M_n(R)$ with
$n\geq2$. If $a_{ii}\geq a_{jk}$ for all $i,j,k\in \underline{n}$,
then \vspace{2mm}

(1)  $\adj(A)=A^{n-1}$;  \vspace{2mm}

 (2)  $\per(\adj(A))=(\per(A))^{n-1}$. }\hfill(3.2) \vspace{3mm}

\noindent {\bf Proof} \quad (1) First, we prove $\per(A(j|i))\leq
a^{(n-1)}_{ij}$ (the $(i,j)$th entry of $A^{n-1}$). If $i=j$, by
assumption, we have
$$\begin{array}{ll}\vspace{2mm} \per(A(i|i))&=\sum\limits_{\sigma\in
S_n}a_{1\sigma(1)}a_{2\sigma(2)}\cdots
a_{i-1,\sigma(i-1)}a_{i+1,\sigma(i+1)}\cdots a_{n\sigma(n)}\\
\vspace{2mm}

 &\leq \sum\limits _{\sigma\in
S_n} a_{11}a_{22}\cdots a_{i-1, i-1}a_{i+1, i+1}\cdots a_{nn}\\
\vspace{2mm}

&=a_{11}a_{22}\cdots a_{i-1, i-1}a_{i+1, i+1}\cdots a_{nn}\\
\vspace{2mm}

&=a_{ii}^{n-1}\\
\vspace{2mm}

 &\leq a^{(n-1)}_{ii} ( \textrm{because } a_{ii}^{n-1} \textrm{ is a term
of } a^{(n-1)}_{ii}).
\end{array}$$
If $i\neq j$, we have
\begin{equation*}
\per(A(j|i))=\sum_{\sigma\in S_n
\atop\sigma(j)=i}a_{1\sigma(1)}a_{2\sigma(2)}\cdots a_{j-1,
\sigma(j-1)}a_{j+1, \sigma(j+1)}\cdots a_{n\sigma(n)}.
\end{equation*}
 Let $T_\sigma =a_{1\sigma(1)}\cdots
a_{j-1,\sigma(j-1)}a_{j+1,\sigma(j+1)}\cdots a_{n\sigma (n)}$ be any
term of $\per(A(j|i))$. Since $\sigma(j)=i$, there exists a positive
integer $u$ such that  $\sigma ^u(i)=j$. Let $t$ be the least
positive integer such that
 $\sigma ^t(i)=j$, then $i, \sigma (i), \ldots , \sigma ^{t-1}(i), j$ are mutually different.
Thus we have $$T_\sigma\leq a_{i\sigma (i)}a_{\sigma (i)\sigma
^2(i)}\cdots a_{\sigma ^{t-1}(i)j}a_{jj}^{n-1-t}\leq
a_{ij}^{(n-1)}.$$ By the choise of  $T_\sigma$, we have
$\per(A(j|i)) \leq a_{ij}^{(n-1)}$, as required.

Next, we prove $a^{(n-1)}_{ij}\leq \per(A(j|i))$. If $i=j$, then we
have
$$\begin{array}{ll}\vspace{2mm} a^{(n-1)}_{ii}&=\sum\limits_{j_1,j_2, \cdots ,
j_{n-2}\in\underline{n}}a_{ij_1}a_{j_1j_2}\cdots a_{j_{n-2}i}\\

 &\leq \sum\limits_{j_1,j_2, \ldots ,
j_{n-2}\in\underline{n}}a_{11}a_{22}\cdots a_{i-1,i-1}a_{i+1,i+1}\cdots a_{nn}\\
\vspace{2mm}

&=a_{11}a_{22}\cdots a_{i-1,i-1}a_{i+1,i+1}\cdots a_{nn}\\

 &\leq \per(A(i|i)),
\end{array}$$
  as required. If $i\neq j$, let $T=a_{ij_1}a_{j_1j_2}\cdots
a_{j_{n-2}j}$ be any term of  the expansion of $a^{(n-1)}_{ij}$. If
$i, j_1, j_2, \ldots, j_{n-2}, j$ are mutually different, then $T$
is a term of $\per(A(j|i))$, which implies $T\leq \per(A(j|i))$. If
$j_r=j_s$ for some $r,s \in \{0, 1, 2,\ldots , n-1 \}$ with $r<s$
(taking $i=j_0$ and $j=j_{n-1}$), delecting $a_{j_rj_{r+1}}\cdots
a_{j_{s-1}j_s}$ from $T$, we get $a_{ij_1}\cdots
a_{j_{r-1}j_r}a_{j_sj_{s+1}}\cdots a_{j_{n-2}j}$. Repeating  the
above method until we can obtain $a_{il_1}a_{l_1l_2}\cdots
a_{l_mj}$, where $0\leq m<n-2$ and $i,l_1,\ldots,l_m,j$ are mutually
different. Let
$\{p_1,p_2,\ldots,p_{n-m-2}\}=\underline{n}\setminus\{i,l_1,\ldots,l_m,j\}$,
we have $$T\leq a_{il_1}a_{l_1l_2}\cdots
a_{l_mj}a_{p_1p_1}a_{p_2p_2}\cdots a_{p_{n-m-2}p_{n-m-2}}\leq
\per(A(j|i)).$$ By the choice of $T$, we have
 $a^{(n-1)}_{ij}\leq\per(A(j|i))$, as required.

 Consequently, $a^{(n-1)}_{ij}=\per(A(j|i))$ which implies
$A^{n-1}=\adj(A)$.

(2)  Using lemma 3.2,  we have $\per(A^{n-1})=(\per A)^{n-1}.$ Thus,
by (1),  we have $\per(\adj(A))=(\per A)^{n-1}$, as required. This
completes the proof.\hfill $\Box$

\vspace{3mm} \noindent {\bf Remark 3.4} \quad Since any commutative
incline is a commutative additively idempotent semiring, Theorem 3.3
solves Problem 1.2.

\vspace{3mm} Wether the equality (3.2) in Theorem 3.3 is hold for
any $n\times n$ matrix over an additively idempotent semiring is not
known. Meanwhile it is well known that for an $n\times n$ matrix $A$
over a field, $|\adj(A)|=|A|^{n-1}$, where $|A|$ is the determinant
of $A$. Recall that  $|A\adj(A)|=|A|^{n}$ for an $n\times n$ matrix
$A$ over a field. For a matrix over a commutative additively
idempotent semiring $R$, we have the following theorem which will be
proved in the following section. \vspace{3mm}

\noindent  {\bf Theorem 3.5} \quad {\itshape For $A\in M_n(R)$ with
$n\geq2$, we have}
$$\per(A\adj(A))=\per(\adj(A)A)=(\per(A))^{n}.\eqno(3.3)$$

 \noindent{\bf Remark 3.6} \quad \emph{Let $R=([0,1],\vee
,\cdot)$, where $\vee=\max$ and $\cdot$ is the usual multiplication
of real numbers. Obviously, $R$ is a commutative incline. Let
$$A=\left(
\begin{array}{ccc}
0.1&0&0.2\\
0&0.2&0.3\\
0&0&0.3
\end{array}
\right)\in M_3(R),$$ we have  $\per(A)=6\cdot10^{-3}$ and
$\per(A\adj(A))=6^{3}\cdot10^{-9}=(\per (A))^3\neq \per(A)$. Since
any commutative incline is a commutative additively idempotent
semiring, Theorem 3.5 solves Problem 1.1 in the negative sense.}

\noindent{\bf Remark 3.7}\quad \emph{Since the fuzzy algebra
$\mathbb{F}=([0,1],\vee,\wedge)$ and any bounded distributive
lattice are commutative
  additively idempotent semiring which satisfying the multiplicatin is idempotent,
  Theorem 3.5 generalizes Theorem 4 in \cite{K} and Theorem 6 in \cite{Z}.}

\vspace{3mm}{\section{{\bf Proof of theorem 3.5 }}}
\setcounter{section}{4} \setcounter{equation}{0}

In the proof, we need the following notation and three lemmas.
\vspace{3mm}

 \noindent {\bf Notation 4.1} \quad (1) For $p,q\in \underline{n},\sigma\in S_{n}$,  denote
 $$\Phi_{p,q}^{\sigma}
=\{(i,\sigma(i))|i\neq
q,i\in\underline{n}\}\cup\{(p,\sigma(q))\},\quad
\Phi^{\sigma}=\Phi_{p,p}^{\sigma}=\{(i,\sigma(i))|i\in\underline{n}\}.$$

(2) For $A\in M_n(R)$, denote by $A(p\Rightarrow q)$  the matrix
obtained from $A$ by replacing row $q$ of $A$ by row $p$ of $A$.
\vspace{3mm}

\noindent {\bf Lemma 4.2} \quad {\itshape Suppose $\sigma,\pi\in
S_{n}$ and $p,q,r\in \underline{n}$ with $q\neq r$, then we have
$$\Phi_{p,q}^{\sigma}\cup\Phi_{q,r}^{\pi}=\Phi^{\varphi}\cup\Phi_{p,r}^{\tau},$$
for some $\varphi,\tau\in S_{n}$. }

  \vspace{3mm}{\bf\noindent Proof.}
  We divide into two cases to define $\varphi$.

\vspace{3mm}\textbf{Case one:} \quad $q \neq(\pi^{-1}\sigma)^k(r)$ for any
nonnegative integer $k$ (we may consider $\sigma ^0$ is the
identical  mapping of the set $\underline{n}$ for any $\sigma \in
S_n$). In this case, define $\varphi:\underline{n}\rightarrow
\underline{n}$ by\
\begin{equation*}\varphi(i)=
\begin{cases}\sigma (i)& \text {if  $i\in U$}\\
             \pi(i) & \text{if $i\in \underline{n}-U$}
\end{cases}
\end{equation*}
\noindent where $U=\{(\pi^{-1}\sigma)^k(r) | k\in \mathbb{N}\}$ and
$\mathbb{N}$ denotes the set of all nonnegative integers. Then,
$\varphi \in S_n$ with $\varphi(q)=\pi(q)$ and
$\varphi(r)=\sigma(r)$. In fact, for any $i,j \in \underline{n}$
with $i\neq j$, if $i,j\in U$ then $\sigma(i)\neq \sigma(j)$ which
implies $\varphi(i)\neq \varphi(j)$, and similarly, if $i,j\in
\underline{n}-U$ then  $\varphi(i)\neq \varphi(j)$. We may suppose
$i\in U$ and $j\in \underline{n}-U$. If $\varphi(i)= \varphi(j)$
then $\sigma(i)=\pi(j)$ and so
 $j=(\pi^{-1}\sigma )(i)$. Since $i=(\pi^{-1}\sigma)^k(r)$ for
 some $k\in \mathbb{N}$, we have $j=(\pi^{-1}\sigma )^{k+1}(r)$
 and so $j\in U$, which is a contradiction. Then $\varphi(i)\neq \varphi(j)$.
Therefore, the mapping $\varphi$ is injective with $\underline{n}$
is a finite set, which implies $\varphi \in S_n$. It is clear that
$\varphi(q)=\pi(q)$ and $\varphi(r)=\sigma(r)$ since $q\in
\underline{n}-U$ and $r\in U$.

 \vspace{3mm}\textbf{Case two:}\quad
$q= (\pi^{-1}\sigma)^k(r)$ for some $k\in \mathbb{N}$. Let $k_0$ be
the least nonnegative integer $k$ such that
$q=(\pi^{-1}\sigma)^k(r)$. Then $q= (\pi^{-1}\sigma)^{k_0}(r)$ and
$k_0\geq 1$ (because $q\neq r$). In this case,  define
$\varphi:\underline{n}\rightarrow \underline{n}$ by \\
\begin{equation*}\varphi(i)=
\begin{cases}\pi(i) & \text {if $i\in U-(V\cup \{q\})$}\\
             \sigma (i)& \text {if  $i\in (\underline{n}-U)\cup
             V$}\\
             \pi(r)& \text {if $i=q$}
\end{cases}
\end{equation*}
\noindent where $V=\{(\pi^{-1}\sigma)^t(r) | t=0,1,\cdots,
k_0-1\}\subseteq U$ (Note that $q\in U-V$). Similarly, we can prove
that $\varphi \in S_n$ with $\varphi (q)=\pi(r)$ and $\varphi
(r)=\sigma (r)$.

\vspace{3mm} By deleting the elements
$(1,\varphi(1)),(2,\varphi(2)),\ldots,(n,\varphi(n))$
from the set $\Phi_{p,q}^{\sigma}\cup\Phi_{q,r}^{\pi}$,
we can get the following $n$ elements:
$$(1,b_1),\ldots,(q,b_q),\ldots,(p,b_r),\ldots,(n,b_n)$$
\noindent where  $b_1, b_2, \ldots,  b_n\in \underline {n}$ are
mutually different, i.e., there exists a permutation $\tau \in S_n$
such that $b_j=\tau(j)$ for all $j\in \underline{n}$. By notation
4.1 (1), we have
$\Phi_{p,q}^{\sigma}\cup\Phi_{q,r}^{\pi}=\Phi^{\varphi}\cup\Phi_{p,r}^{\tau}$
as required. This completes the proof. \hfill $\Box$

 \vspace{3mm}\noindent {\bf Lemma 4.3}\quad
{\itshape Let $A\in M_n(R)$ and $k\in \underline{n}$ with $n\geq2$.
We have
$$\per(A(p_1\Rightarrow p_2))\per(A(p_2\Rightarrow p_3))\cdots \per(A(p_k\Rightarrow
p_1))\leq(\per(A))^{k},$$ for $p_1,p_2,\ldots,p_k\in
\underline{n}$.}

 \vspace{3mm}{\bf\noindent Proof.}  We first prove that for $p,q,r \in
\underline{n}$,
$$\per(A(p\Rightarrow q))\cdot\per(A(q\Rightarrow
r))\leq\per(A)\cdot\per(A(p\Rightarrow r)).\eqno(4.1)$$

 In fact, the case $q=r$ is obvious by commutativity. Suppose $q\neq r$, then
$$\begin{array}{ll}\vspace{2mm} &\per(A(p\Rightarrow q))\cdot \per(A(q\Rightarrow
r))\\ \vspace{2mm}

 =&\sum\limits_{\sigma \in S_n}a_{1\sigma (1)}\cdots
a_{p\sigma (q)}\cdots a_{n\sigma (n)}\cdot \sum\limits_{\pi \in
S_n}a_{1\pi(1)} \cdots a_{q\pi(r)}\cdots a_{n\pi(n)}\\
 \vspace{2mm}

 =&\sum\limits_{\sigma,\pi \in S_n}a_{1\sigma (1)}\cdots a_{p\sigma
(q)}\cdots a_{n\sigma (n)}a_{1\pi(1)}\cdots a_{q\pi(r)}\cdots
a_{n\pi(n)}.\end{array}$$ Denote by  $(i,j)$ the subscript of
$a_{ij}$, we get the set $\Phi_{p,q}^{\sigma}\cup\Phi_{q,r}^{\pi}$
for $$T:=a_{1\sigma (1)}\cdots a_{p\sigma (q)}\cdots a_{n\sigma
(n)}a_{1\pi(1)}\cdots a_{q\pi(r)}\cdots a_{n\pi(n)}.$$ By lemma 4.2
we have
$\Phi_{p,q}^{\sigma}\cup\Phi_{q,r}^{\pi}=\Phi^{\varphi}\cup\Phi_{p,r}^{\tau}$
 for some $\varphi,\tau\in S_{n}$. Therefore
$$T=a_{1\varphi(1)}a_{2\varphi(2)}\cdots a_{n\varphi (n)}a_{1\tau(1)}\cdots a_{p\tau(r)}\cdots
a_{n\tau(n)}$$ is a term of $\per(A)\per(A(p\Rightarrow r))$, which
implies $ T \leq\per(A) \cdot \per(A(p\Rightarrow r)).$ Since $T$ is
any term of $\per(A(p\Rightarrow q))\cdot \per(A(q\Rightarrow r))$,
we have (4.1).

Using (4.1), we have
$$\begin{array}{ll}\vspace{2mm}
&\per(A(p_1\Rightarrow p_2))\per(A(p_2\Rightarrow p_3))\cdots
\per(A(p_k\Rightarrow p_1))\\ \vspace{2mm}

\leq&\per(A)\cdot \per(A(p_1\Rightarrow p_3))\per(A(p_3\Rightarrow
p_4))\cdots \per(A(p_k\Rightarrow p_1))\\
\vspace{2mm}

 \leq&(\per(A))^2\cdot \per(A(p_1\Rightarrow
p_4))\per(A(p_4\Rightarrow
p_5))\cdots \per(A(p_k\Rightarrow p_1))\\
\vspace{2mm}

&\cdots\cdots\cdots\cdots\\
\vspace{2mm}

\leq&(\per(A))^{(k-1)}\cdot \per(A(p_1\Rightarrow p_1))\\
=&(\per(A))^k,
\end{array}$$
as required. This completes the proof.\hfill $\Box$ \vspace{3mm}

\noindent {\bf Lemma 4.4}\quad {\itshape For $A\in M_n(R)$ with
$n\geq2$, we have}
$$\per(A\adj(A))\geq(\per(A))^{n}.$$

\noindent  {\bf Proof. }   Denote by $b_{ij}$ the $(i,j)$-th element
of $A\adj(A)$, we have
$$b_{ij}=\sum\limits_{k=1}^n a_{ik}\per(A(j|k))=\per(A(i\Rightarrow j)).$$
Thus we have
$$
\begin{array}{ll}\vspace{2mm}

 &\per(A\adj(A))\\
 \vspace{2mm}

=&\sum\limits _{\pi\in S_n}\per(A(1\Rightarrow
\pi(1)))\per(A(2\Rightarrow \pi(2)))\cdots \per(A(n\Rightarrow
\pi(n)))\\
\vspace{2mm}

\geq&\per(A(1\Rightarrow 1))\per(A(2\Rightarrow 2))\cdots
\per(A(n\Rightarrow
n))\\
\vspace{2mm}

=&(\per(A))^n,
\end{array}$$
as required. This completes the proof. \hfill $\Box$

\vspace{3mm}\noindent  \textbf{Proof of Theorem 3.5} \quad By lemma
4.4, we have $\per(A\adj(A))\geq(\per(A))^{n}$. In the following, we
will prove $\per(A\adj(A))\leq(\per(A))^{n}$.

Let $T:=\per(A(1\Rightarrow \pi(1)))\per(A(2\Rightarrow
\pi(2)))\cdots \per(A(n\Rightarrow \pi(n)))$ be any term of
$\per(A\adj(A))$, where $\pi\in S_n$. Let $\pi =(i_1i_2\cdots
i_r)(j_1j_2\cdots j_s)\cdots (l_1l_2\cdots l_u)$ be the
decomposition of disjoint cycles with $r+s+\cdots +u=n$. Using lemma
4.3, we have
$$\begin{array}{ll}\vspace{2mm}

T&=\prod\limits_{k\in\underline{r}}\per(A(i_k\Rightarrow
\pi(i_k)))\prod\limits_{k\in\underline{s}} \per(A(j_k\Rightarrow
\pi(j_k)))\cdots
\prod\limits_{k\in\underline{u}}\per(A(l_k\Rightarrow \pi(l_k)))\\
\vspace{2mm}
 &=\big(\prod\limits_{ k
\in\underline{r-1}}\per(A(i_k\Rightarrow
i_{k+1}))\big)\per(A(i_r\Rightarrow i_1))\cdot
\big(\prod\limits_{ k \in\underline{s-1}}\per(A(j_k\Rightarrow j_{k+1}))\big)\\
\vspace{2mm}

&\quad \cdot\per(A(j_s\Rightarrow j_1)) \cdots  \big(\prod\limits_{
k \in\underline{u-1}}\per(A(l_k\Rightarrow
l_{k+1}))\big) \per(A(l_u\Rightarrow l_1))\\
\vspace{2mm}

&\leq (\per(A))^{r}(\per(A))^{s}\cdots (\per(A))^{u}\\
\vspace{2mm}

&=(\per(A))^{r+s+\cdots+u}=(\per(A))^{n}.\end{array}$$ By the choice
of $T$, we have $\per(A\adj(A))\leq (\per A)^{n}$ as required.
Therefore, we have $\per(A\adj(A))=(\per(A))^{n}$.

Finally, using proposition 2.1(2) and proposition 3.1(3),  we have
$\per(\adj(A)A)=(\per(A))^{n}.$ This completes the proof.\hfill
$\Box$


\begin{thebibliography}{99}
\setlength{\baselineskip}{1.1em}
\bibitem{B} J.P.M. Binet, M\'{e}moire sur un syst\`{e}me de formulas analytiques, et
leur application \`{a} des consid\'{e}rations g\'{e}ometriques, J.
\'{E}c. Polyt., 9 (1812), 280-302.
\bibitem{Du} J.S. Duan, The transitive closure, convergence of powers and adjoint of generalized fuzzy
matrices,
  Fuzzy Sets and System, 145 (2004), 301-311.
\bibitem{G}D.G. Glynn, The permanent of a square matrix, European Journal of Combinatorics, 31 (2010), 1887-1891.
\bibitem{GM} M. Gondran, M. Minoux, Dio\"{\i}ds and semirings: links to fuzzy sets and other
applications, Fuzzy Sets and Systems, 158 (2007) 1273-1294.
\bibitem{H1} S. Han, H. Li, Invertible incline matrices and Cramer's rule
over inclines, Linear Algebra and its Applications, 389 (2004)
121-138.
\bibitem{K} J.B. Kim, A. Baartmans, N.S. sahadin, Determinant theory for fuzzy matrices, Fuzzy Sets and System, 29 (1989), 349-356.
\bibitem{KM} E.P. Klement, R. Mesiar, E. Pap, Triangular Norms, Kluwer, Dordrecht, 2000.
\bibitem{LB} H. Liang, F. Bai, An upper bound for the permanent of (0,1)-matrices, Linear Algebra and its Applications, 377 (2004), 291-295.
\bibitem{M} H. Minc, Permanents, Addison-Wesley Publishing Company, Massachusetts, U. S.
A., 1978.
 \bibitem{RE} M.Z. Ragab, E.G. Emam, The determinant and adjoint of a square fuzzy matrix, Fuzzy Sets and Systems, 61 (1994), 297-307.
 \bibitem{S} G. Sburlati, On the parity of permanents of circulantmatrices, Linear Algebra and its Applications,  428 (2008),
1949-1955.
\bibitem{Th} M.G. Thomason, Convergence of powers of a fuzzy matrix, J.
Math. Anal. Appl. 57 (1977) 476¨C480.
\bibitem{Z} K.L. Zhang, Determinant theory for $D_{01}$-lattice matrices, Fuzzy Sets Systems, 62 (1994), 347-353.


 \end{thebibliography}
\end{document}